%% file: main.tex
\documentclass[preprint,12pt]{elsarticle}

\usepackage{amssymb}
\usepackage{amsmath}

\usepackage{nomencl}
\makenomenclature

\usepackage{etoolbox}
\renewcommand\nomgroup[1]{%
  \item[\bfseries
  \ifstrequal{#1}{S}{Superscripts}{%
  \ifstrequal{#1}{P}{Parameters}{%
  \ifstrequal{#1}{I}{Sets and Indices}{%
  \ifstrequal{#1}{V}{Variables}{%
  \ifstrequal{#1}{F}{Functions}{}}}}}%
]}

\usepackage{xcolor}

\usepackage{multirow}
\usepackage{hyperref}
\usepackage[acronym]{glossaries}

\journal{Sustainable Energy, Grids and Network}

\begin{document}

\begin{frontmatter}

\title{Understanding the Impact of Hydro-Reservoirs and Inverters on Frequency-Constrained Operation}

\author[UC]{Valeria Aravena}
\author[UC]{Samuel Córdova\corref{cor1}}
\cortext[cor1]{Corresponding author}
\ead{sacordov@uc.cl}
\author[UC]{Maximiliano Kairath}
\author[UC]{Matías Negrete-Pincetic}

\affiliation[UC]{organization={Pontifical Catholic University of Chile},
            addressline={Vicuña Mackenna 4860}, 
            city={Santiago},
            postcode={7820436}, 
            country={Chile}}

\begin{abstract}
The increasing participation of renewable energy sources in power systems has entailed a series of challenges resulting from the replacement of conventional synchronous machines with carbon-free Inverter-Based Resources (IBRs). In this context, the present work contributes to the existing literature on Frequency-Constrained Unit Commitment (FCUC) models by studying the role of hydro-reservoirs in the ongoing decarbonization of power systems.  For this purpose, a novel FCUC model is developed, which captures hydro-reservoir dynamics and their impact on frequency Nadir requirements through a data-driven approach. Moreover, the proposed FCUC model is used to simulate a series of future scenarios in the Chilean power system, and to understand the role that hydropower units, thermal units, Grid-Forming (GFM) inverters, and Synchronous Condensers (SCs) will play in the future in terms of frequency regulation.
Exhaustive simulations on the Chilean power system for years 2024 (current scenario) and 2035 (carbon-free scenario) illustrate the benefits of the proposed approach, which include: (i) significant cost-savings resulting from using the proposed FCUC model relative to the industry standard, achieving operational cost reductions of up to 28\% without compromising system security; and 
(ii) providing insights into the key role that hydro-reservoirs will play in the future.
As a final analysis, the developed FCUC model is used to evaluate the financial benefits of incorporating SCs and GFM inverters to provide frequency regulation support to the grid in 2035. The results show that, while both technologies contribute to reducing the system's annual operational costs, GFM inverters significantly outperform SCs in terms of investment return.
\end{abstract}



\end{frontmatter}





\section{Introduction}
\label{sec:intro}
\input{./introduction.tex}

\section{Methodology}
\label{sec:methodology}
\input{./methodology.tex}

\section{Computational Experiments}
\label{sec:computational_experiments}
\input{./computational_experiments.tex}

\section{Conclusion}
\label{sec:conclusion}
\input{./conclusion.tex}


\section*{CRediT authorship contribution statement}
\input{./CRediT.tex}

\bibliographystyle{IEEEtran} 
\bibliography{references}

\end{document}

%% file: introduction.tex
The ongoing decarbonization of power systems around the globe has entailed the displacement of conventional thermal-based Synchronous Machines (SMs) by modern carbon-free Inverter-Based Resources (IBRs), such as wind, solar PV, and battery storage plants. However, this transition introduces several challenges, most notably a reduction in available system inertia. As a result, frequency regulation during contingencies is diminished \cite{Denholm2020}, increasing the risk of grid instability \cite{Badesa2019}. 

This new generation mix has led to the need for novel Frequency-Constrained Unit Commitment (FCUC) models suited to IBR-dominated system operation \cite{Paturet2020}. This ensures that grid frequency remains within predefined frequency limits, which are associated with the metrics of Rate of Change of Frequency (RoCoF), Nadir, and Quasi-Steady State (QSS). 
Unlike RoCoF and QSS, the Nadir metric is characterized by its non-linear dependence on various operational variables \cite{Paturet2020, Markovic}. To address this challenge, the literature often employs analytical methodologies, which are generally based on linearizations (e.g., \cite{Ahmadi2014, Teng2016, Shi2025, Zeng2023}). However, the approximation error associated with this approach can result in suboptimal or even infeasible solutions. As an alternative, data-driven methodologies have been proposed, in which Nadir constraints are inferred through exhaustive dynamic simulations (e.g., \cite{Shen2023, Rajabdorri2023, Sebastian2024, Wu2024, ZhangC2021, Zhang2023}). 

\begin{figure}[h]
    \centering
    \includegraphics[width=0.6\linewidth]{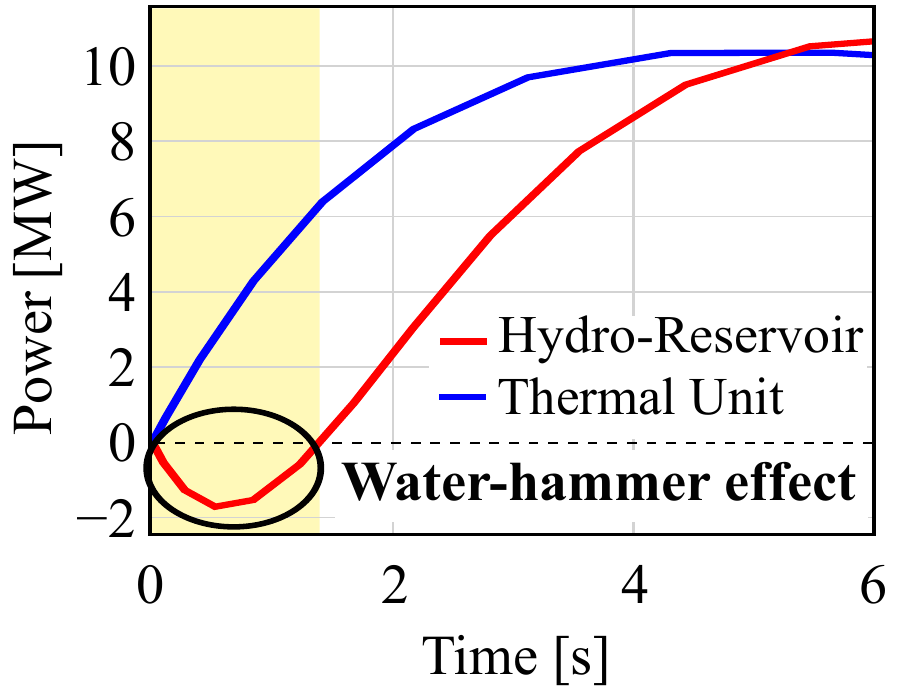}
    \caption{Power response of a thermal and a hydro-reservoir unit}
    \label{fig:Figure_1}
\end{figure}

While significant research efforts have been put into developing FCUC models that integrate both thermal units and IBRs (e.g., \cite{Zhou2024, Zhang2020, Paturet2020, Wogrin2020, Shen2023, Rajabdorri2023, Ferrandon-Cervantes2022, He2024}), less attention has been given to understanding the impact of hydro-reservoirs and their dynamics on the problem (e.g., \cite{Sebastian2024, ZhangQ2021, ZhangC2021, Zhang2023}). This can be particularly relevant for two reasons: (i) countries like China, Chile, and Brazil have a significant share of hydro-reservoirs in their total generation mix (e.g., 16\% of total generation in Chile during 2024 \cite{CEN2025}), making their dynamic response systemically relevant; and (ii) hydro-reservoir generation plants are affected by the water-hammer effect, a hydraulic phenomenon that causes the power output to initially respond in the opposite direction to that desired \cite{Kundur}, which can impact compliance with frequency metrics (see Fig.~\ref{fig:Figure_1}). 
Moreover, state-of-the-art FCUC models now typically include modern technologies for frequency support, such as Grid-Forming (GFM) inverters (e.g. \cite{Paturet2020, Sebastian2024}).

Reference \cite{ZhangQ2021} presents an FCUC model that considers the operation of hydro-reservoirs and GFM inverters, in which the Nadir constraint is incorporated using an analytical methodology. This constraint, however, is obtained by linearizing an expression derived from a thermal system model, thereby omitting the effect of the dynamic response of hydro-reservoirs. On the other hand, reference \cite{Sebastian2024} employs a data-driven methodology, where simulations are performed for different levels of system inertia and droop. While hydro-reservoir operation is considered, its dynamic response is not included in the model, underestimating its impact on frequency regulation. Similarly, \cite{ZhangC2021} also follows a data-driven approach and incorporates hydro-reservoirs in its analysis. Nevertheless, the water-hammer effect is neglected, resulting in an inaccurate representation of hydro-reservoir dynamics.
The above problem is addressed in \cite{Zhang2023}, which uses a complete dynamic representation of hydro-reservoir plants including the water-hammer effect. In this case, the Nadir constraint is also formulated using a data-driven approach, with simulations performed for different levels of inertia, contingency sizes, and online capacity of generation technologies. Based on the results, the conditions under which the constraint is satisfied or violated are identified, and a Support Vector Classification algorithm is used to generate cutting planes that separate the two regions of constraint satisfaction. These planes are then integrated as linear constraints into the FCUC model, which allows for accurately capturing the dynamics of hydro-reservoirs. Nevertheless, the above work neglects the potential of GFM inverters for frequency regulation provision. 

Based on the above literature review, a research gap is identified in the development of an FCUC model that integrates the dynamic response of hydro-reservoirs and GFM inverters. Accordingly, the main objectives and contributions of this paper are:

\begin{enumerate}
    \item Simulate the dynamic response of different generation and/or grid-supporting technologies (e.g., GFM inverters) and evaluate the compliance of RoCoF, Nadir, and QSS frequency requirements.
    \item Develop an FCUC model that incorporates the dynamic response of hydro-reservoirs and GFM inverters, and evaluate its application for the Chilean electrical system, characterized by a significant participation of solar and hydro resources.
    \item Perform exhaustive computational simulations to understand the role of conventional thermal/hydro generators and emerging GFM technologies in frequency regulation for the next years in Chile. Moreover, based on such simulations, propose general guidelines for future carbon-free system operation. 
\end{enumerate}

The rest of the paper is organized as follows: Section \ref{sec:methodology} presents a frequency dynamic model incorporating hydro-reservoir and GFM, which is then used to develop a novel FCUC model embedding such dynamics using a data-driven approach. Next, Section \ref{sec:computational_experiments} describes the computational experiments, detailing the results related to the system’s dynamic response, the evaluation of the developed model, and the integration of GFM inverters and Synchronous Condensers (SCs) into the future operation of the Chilean electrical system. Finally, Section \ref{sec:conclusion} presents the main findings of the research.

%% file: methodology.tex
This section presents the proposed dynamic and FCUC models used for the analysis. On the one hand, Section \ref{sec:dynamic_model} describes the frequency dynamic model, which includes the dynamic response of conventional generators, the water-hammer effect of hydro-reservoirs, and GFM capabilities for inverter-based resources. On the other hand, Section \ref{sec:FCUC} presents the novel FCUC model, which includes the dynamics of hydro-reservoirs and GFM inverters.

\subsection{Dynamic Modeling of Power Systems} \label{sec:dynamic_model} 

\begin{figure}[h] 
    \centering
    \includegraphics[width=1\linewidth]{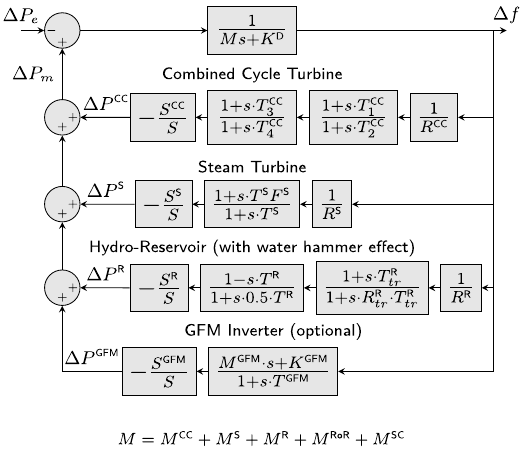}
    \caption{System frequency dynamics model}
    \label{fig:Figure_2}
\end{figure}

Figure \ref{fig:Figure_2} presents a dynamic model of a system in which steam turbines (coal-fired power plants), combined-cycle turbines (gas-fired power plants), hydro-reservoirs, and GFM inverters participate in primary frequency control. In the model, $\Delta P_e$ denotes the contingency size, $M$ the system inertia, $\Delta f$ the frequency variation, $S$ the system base power, while $S^\text{CC}$, $S^\text{S}$, $S^\text{R}$, $S^\text{GFM}$ denote the online capacity of combined-cycle turbine plants, steam turbine plants, hydro-reservoirs, and GFM inverters, respectively. The model also accounts for load damping $K^D$, as well as the inertial response of generation technologies participating in primary control ($M^\text{CC}$, $M^\text{S}$, $M^\text{R}$, $M^\text{GFM}$), run-of-river hydropower plants ($M^\text{RoR}$), and synchronous condensers ($M^\text{SC}$). 

Since multiple technologies participate in the system, an aggregated representation known as the Center of Inertia (CoI) is used, where $M$ is the weighted sum of the inertia provided by the different units. Regarding the primary control models for the different technologies, a first-order combined cycle turbine model is employed, adapted from the one presented in \cite{FERNANDEZGUILLAMON}. The models for steam turbines and hydro-reservoirs are adopted from \cite{Kundur}. For GFM inverters, the model is based on \cite{Paturet2020} and corresponds to a Virtual Synchronous Machine (VSM) configuration.

\subsection{Frequency-Constrained Unit Commitment} \label{sec:FCUC}
This section presents the objective function and the constraints associated with a daily operation model that considers both standard constraints (Section \ref{sec:FCUC_standar}) and frequency constraints (Section \ref{sec:FCUC_frequency}). 

In this model, $t \in \mathcal{T}$, $i \in \mathcal{I}$, $h \in \mathcal{H}$, $b \in \mathcal{B}$, and $r \in \mathcal{R}$ represent the indices and sets for time periods, thermal units, hydropower units, batteries, and renewable units, respectively. Within hydropower units, two subsets are distinguished: reservoir-based hydropower plants ($\mathcal{H}_R$) and run-of-river hydropower plants ($\mathcal{H}_{RoR}$). Likewise, batteries are classified into those equipped with Grid-Following (GFL) inverters ($\mathcal{B}_{GFL}$) and those equipped with GFM inverters ($\mathcal{B}_{GFM}$). 

The objective function is defined as follows:
\begin{equation}
    \label{ec:FO}
    \begin{split}
        \min & \sum_{i \in \mathcal{I}} \sum_{t \in \mathcal{T}} \left( C_{i}^{V} p_{it} + C_{i}^{F} u_{it} + C_{i}^{SU} y_{it} + C_{i}^{SD} z_{it} \right) \\
        & + \sum_{h \in \mathcal{H}_E} \sum_{t \in \mathcal{T}} C_{h}^{V} p_{ht} \\
         & + \sum_{b \in \mathcal{B}} \sum_{t \in \mathcal{T}} \left( C_b^{V} p^{dis}_{bt} + 
         C_b^{V} p^{cha}_{bt} \right)
    \end{split}
\end{equation}
where the costs associated with thermal operation include: variable costs $C^V_{it}$, which depend on the output power of each unit $p_{it}$; fixed costs $C^F_{it}$, determined by the binary commitment variable $u_{it}$; start-up costs $C^{SU}_{it}$, associated with the binary start-up variable $y_{it}$; and shut-down costs $C^{SD}_{it}$, associated with the binary shut-down variable $z_{it}$.
For reservoir hydropower, a variable cost $C^{V}_{ht}$ is considered, which is determined by the output power $p_{ht}$ and represents the opportunity cost of using the water resource now versus storing it for future use \cite{valor_agua_CEN}. In the case of batteries, the variable cost $C^{V}_{bt}$ is a function of the charging $p^{cha}_{bt}$ and discharging $p^{dis}_{bt}$ power, and is linked to unit degradation resulting from charge and discharge cycles, reflecting the reduction in useful life due to intensive operation of the storage system \cite{Córdova}.

\subsubsection{Standard constraints} \label{sec:FCUC_standar}

The constraint to maintain the system’s energy balance is expressed as follows:
\begin{equation}
    \label{ec:balance}
    \sum_{g \in \mathcal{G}} p_{gt} + \sum_{b \in \mathcal{B}} (p^{dis}_{bt} - p^{cha}_{bt}) = D_t \quad \forall t
\end{equation}
where the total generation must meet the demand $D_t$ at each time period $t$. \\

Then, the constraints referring to the operation of thermal and hydro-reservoir units are represented as follows:
\begin{equation}
    \label{ec:conv_1}
    y_{gt} - z_{gt} = u_{gt} - u_{g,t-1} \quad \forall g \in (\mathcal{I} \cup \mathcal{H}_R),t\geq2 
\end{equation}
\begin{equation}
    \label{ec:conv_2}
    p_{gt} + r^{+}_{gt} \leq \overline{P}_g u_{gt} \quad \forall g \in (\mathcal{I} \cup \mathcal{H}_R),t
\end{equation}
\begin{equation}
    \label{ec:conv_3}
    p_{gt} \geq \underline{P}_g u_{gt} \quad \forall g \in (\mathcal{I} \cup \mathcal{H}_R),t
\end{equation}
 \begin{equation}
\label{ec:conv_4}
    \sum_{t}^{t+\underline{T}^{ON}_i} u_{it} \geq \underline{T}^{ON}_i y_{it} \quad \forall i ,t \leq 24 -\underline{T}^{ON}_i
\end{equation}
\begin{equation}
\label{ec:conv_5}
    \sum_{t}^{t+\underline{T}^{OFF}_i-1} (1 - u_{it}) \geq \underline{T}^{OFF}_i z_{it} \quad \forall i,t \leq 24 -\underline{T}^{OFF}_i
\end{equation}
\begin{equation}
\label{ec:conv_6}
    \sum_{t}^{24} u_{it} \geq (24-t+1) y_{it} \quad \forall i,t > 24 -\underline{T}^{ON}_i
\end{equation}
\begin{equation}
\label{ec:conv_7}
    \sum_{t}^{24} (1-u_{it}) \geq (24-t+1) z_{it} \quad \forall i,t > 24 -\underline{T}^{OFF}_i
\end{equation}
\begin{equation}
\label{ec:conv_8}
    \sum_{t} p_{ht} \leq E_h \quad \forall h \in \mathcal{H}_R
\end{equation}
where (\ref{ec:conv_1}) defines the logic for the binary variables of commitment $u_{it}$, start-up $y_{it}$, and shut-down $z_{it}$; (\ref{ec:conv_2})-(\ref{ec:conv_3}) set the upper $\overline{P}_g$ and lower $\underline{P}_g$ bounds of power output, accounting for the reserve commitments $r^{+}_{gt}$ of the units; (\ref{ec:conv_4})–(\ref{ec:conv_7}) model the minimum up-time $\underline{T}^{ON}_i$ and down-time $\underline{T}^{OFF}_i$ requirements for thermal units; and (\ref{ec:conv_8}) limits the daily energy generation in $E_h$ for each hydro-reservoir unit. 

For renewable and run-of-river hydropower units, the operational constraints are formulated as follows:
\begin{equation}
\label{ec:renew_1}
    p_{gt} \leq PF_{gt}\ \quad \forall g \in (\mathcal{R} \cup \mathcal{H}_{RoR}),t
\end{equation}
\begin{equation}
\label{ec:renew_2}
    p_{gt} \geq \underline{P}_{g} \quad \forall g \in (\mathcal{R} \cup \mathcal{H}_{RoR}),t
\end{equation}
where (\ref{ec:renew_1}) describes the maximum available power over time $PF_{gt}$; and (\ref{ec:renew_2}) sets the lower bound of power output in $\underline{P}_{g}$. 

The operation of batteries is described by the following constraints:
\begin{equation}
\label{ec:bess_1}
    p^{dis}_{bt} + r^{dis,+}_{bt} \leq \overline{P}_{b} \quad\forall  b,t
\end{equation}
\begin{equation}
\label{ec:bess2}
    r^{cha,+}_{bt} \leq p^{cha}_{bt} \quad \forall  b,t
\end{equation}
\begin{equation}
\label{ec:bess_3}
    p^{cha}_{bt} \leq \overline{P}_{b} \quad\forall  b,t
\end{equation}
\begin{equation}
\label{ec:bess_4}
    p^{dis}_{bt} \geq 0, p^{cha}_{bt} \geq 0 \quad \forall  b,t
\end{equation} 
\begin{equation}
\label{ec:bess_5}
    r^{cha,+}_{bt} = r^{dis,+}_{bt} = 0 \quad \forall  b \in \mathcal{B}_{GFL},t
\end{equation} 
\begin{equation}
\label{ec:bess_6}
    e_{bt} - e^{ini}_{b} = \left( p^{cha}_{bt} \eta^{cha}_{b} - \frac{p^{dis}_{bt}}{\eta^{dis}_{b}} \right) \Delta t \quad\forall   b,t=1
\end{equation}
\begin{equation}
\label{ec:bess_7}
    e_{bt} - e_{b,t-1} = \left( p^{cha}_{bt} \eta^{cha}_{b} - \frac{p^{dis}_{bt}}{\eta^{dis}_{b}} \right) \Delta t \quad\forall  b,t
\end{equation}
\begin{equation}
\label{ec:bess_8}
    e_{bt} \leq \overline{E}_b \quad\forall b,t
\end{equation}
\begin{equation}
\label{ec:bess_9}
    e_{bt} \geq 0 \quad\forall b,t
\end{equation}
\begin{equation}
\label{ec:bess_10}
    e_{bt} - \left( r^{cha,+}_{bt} \eta^{cha}_{b} + \frac{r^{dis,+}_{bt}}{\eta^{dis}_{b}} \right) \Delta t \geq \underline{E}_b \quad\forall b,t
\end{equation}
\begin{equation}
\label{ec:bess_11}
    e_{b,24} = e_{b}^{init} \quad\forall b
\end{equation}
where (\ref{ec:bess_1})–(\ref{ec:bess_4}) define the upper $\overline{P}_{b}$ and lower $\underline{P}_{b}$ bounds of charging and discharging power while accounting for upward reserve provision following the approach in \cite{Córdova}, with storage units providing upward reserves in two ways, either by discharging upward reserve $r^{dis,+}_{bt}$ which involves temporarily increasing the discharge power or by charging upward reserve $r^{cha,+}_{bt}$ which consists of reducing the charging power; (\ref{ec:bess_5}) enforces that reserve provision is zero for batteries operating with GFL inverters, reflecting the assumption that only units equipped with GFM inverters are capable of providing reserves; (\ref{ec:bess_6})–(\ref{ec:bess_7}) model the stored energy of the storage units at each hour $e_{bt}$, considering charging $\eta^{cha}_{b}$ and discharging $\eta^{dis}_{b}$ efficiency; (\ref{ec:bess_8})–(\ref{ec:bess_10}) set the upper $\overline{E}_b$ and lower $\underline{E}_b$ bounds on this stored energy; finally, (\ref{ec:bess_11}) ensures that the energy level at the end of the day matches the one had at the beginning of the day $e^{init}_b$ (initial value).

Furthermore, constraints referring to the system's reserve requirement are defined as follows \cite{Córdova}:
\begin{equation}
\label{ec:reserves_1}
    \sum_{s \in \mathcal{S}} r^{+}_{st} + \sum_{b \in \mathcal{B}_{GFM}} (r^{dis,+}_{bt} + r^{cha,+}_{bt}) + r^{PF+}_t \geq \Delta P_e \quad \forall t
\end{equation}
\begin{equation}
\label{ec:reserves_2}
    \begin{split}
        & r^{+}_{st} \geq 0, r^{dis,+}_{bt} \geq 0, r^{cha,+}_{bt} \geq 0, r^{PF+}_t \geq 0 \quad \\ 
        & \forall s,b \in \mathcal{B}_{GFM},t
    \end{split}
\end{equation}
where (\ref{ec:reserves_1}) guarantees that all upward reserves are enough to cover the power imbalance $\Delta P_e$, accounting for the reserve commitments of generation units as well as the reserve associated with load damping ($r^{PF+}_t$); and (\ref{ec:reserves_2}) requires that all reserves be non-negative.

\subsubsection{Frequency Constraints} \label{sec:FCUC_frequency} 
\paragraph {i} Rate of Change of Frequency (RoCoF) 

According to \textbf{\cite{Kundur}}, the constraint that ensures RoCoF meets the minimum threshold $\dot{f}_{lim}$ is as follows:
\begin{equation}
\label{ec:rest_rocof}
    \frac{\Delta P_e}{m_t} \leq \frac{\dot{f}_{lim}}{f_0} \quad \forall t
\end{equation}
where $m_t$ denotes the system inertia and $f_0$ is the nominal frequency. The total system inertia $m_t$ is determined by the following constraint:
\begin{equation}
\begin{split}
\label{ec:rest_inercia}
      m_t = & \sum_{i \in \mathcal{I}} 2 \cdot H_i \cdot \overline{P}_i \cdot u_{it}  + \sum_{h \in \mathcal{H}_R} 2 \cdot H_h \cdot \overline{P}_h \cdot u_{ht} \\
      & + \sum_{h \in \mathcal{H}_{RoR}} 2 \cdot H_h \cdot \overline{P}_h + \sum_{b \in \mathcal{B}_{GFM}} 2 \cdot H_b \cdot \overline{P}_b \\
      & + \sum_{i \in \mathcal{CS}} 2 \cdot H_i \cdot \overline{P}_i 
      \quad \forall t
\end{split}
\end{equation}
which accounts for the contribution of dispatched thermal and hydro-reservoir units, as well as that from run-of-river hydropower plants, GFM VSM inverters, and SCs available in the grid.
 
\paragraph {ii} Quasi-Steady State (QSS) 

The constraints associated with QSS are based on those presented in \cite{Córdova} and are as follows:
\begin{equation}
\label{ec:qss_1}
    r^{+}_{it} \leq \frac{f_{ss,lim}}{f_0} \frac{\overline{P}_i}{R_i} \quad \forall i,t
\end{equation}
\begin{equation}
\label{ec:qss_2}
    r^{+}_{ht} \leq \frac{f_{ss,lim}}{f_0} \frac{\overline{P}_h}{R_h} \quad \forall h \in \mathcal{H}_E,t
\end{equation}
\begin{equation}
\label{ec:qss_3}
    r^{dis,+}_{bt} + r^{cha,+}_{bt} \leq \frac{f_{ss,lim}}{f_0} \frac{\overline{P}_b}{R_b} \quad \forall b \in \mathcal{B}_{GFM},t
\end{equation}
\begin{equation}
\label{ec:qss_4}
    r^{PF+}_t \leq \frac{f_{ss,lim}}{f_0} K^{D} \quad \forall t
\end{equation}
where frequency-dynamic limits are imposed on the reserves to ensure that the QSS complies with the minimum threshold ${f}_{ss,lim}$. For generating units (\ref{ec:qss_1})-(\ref{ec:qss_3}), the reserve limit depends on the droop constant $R_g$ characteristic of each unit. In the case of demand (\ref{ec:qss_4}), its contribution is determined by the load damping constant $K^D$.

\paragraph {iii} Nadir 

Unlike the previous metrics, the Nadir expression depends on a nonlinear relationship among various system factors, including the magnitude of the contingency $\Delta P_e$, the load damping constant $K^D$, and the system inertia $m_t$ \cite{Kundur}. To manage this complexity, the Nadir constraint is integrated using a data-driven approach.

\begin{figure}[h]
    \centering
    \includegraphics[width=0.65\linewidth]{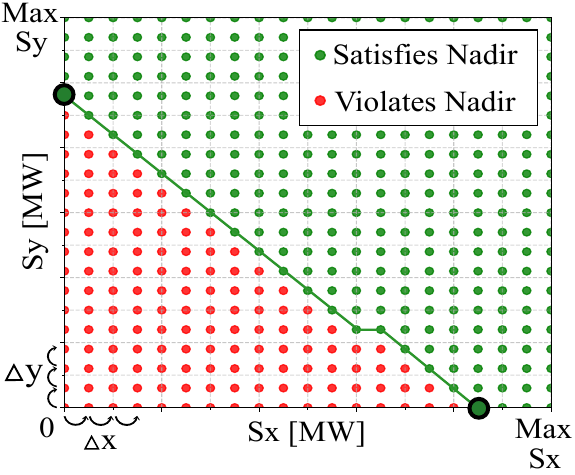}
    \caption{Method used for determining the Nadir Boundary based on generation technologies X and Y. Variables $S_x$ and $S_y$ denote the online capacities of technologies X and Y, respectively. $\Delta x$ and $\Delta y$ indicate the granularity used in the dynamic simulations for determining the Nadir Boundary.}
    \label{fig:Figure_3}
\end{figure}

A sweep of simulations using the dynamic model presented in Figure \ref{fig:Figure_3} is performed, considering different combinations of online capacity for the technologies. For each combination, it is determined whether the Nadir requirement is satisfied, which is indicated in Figure \ref{fig:Figure_3} through green and red points, respectively. Based on these simulations, the Nadir compliance boundary is determined, which is indicated through the solid green line.

\begin{figure}[h]
    \centering
    \includegraphics[width=0.8\linewidth]{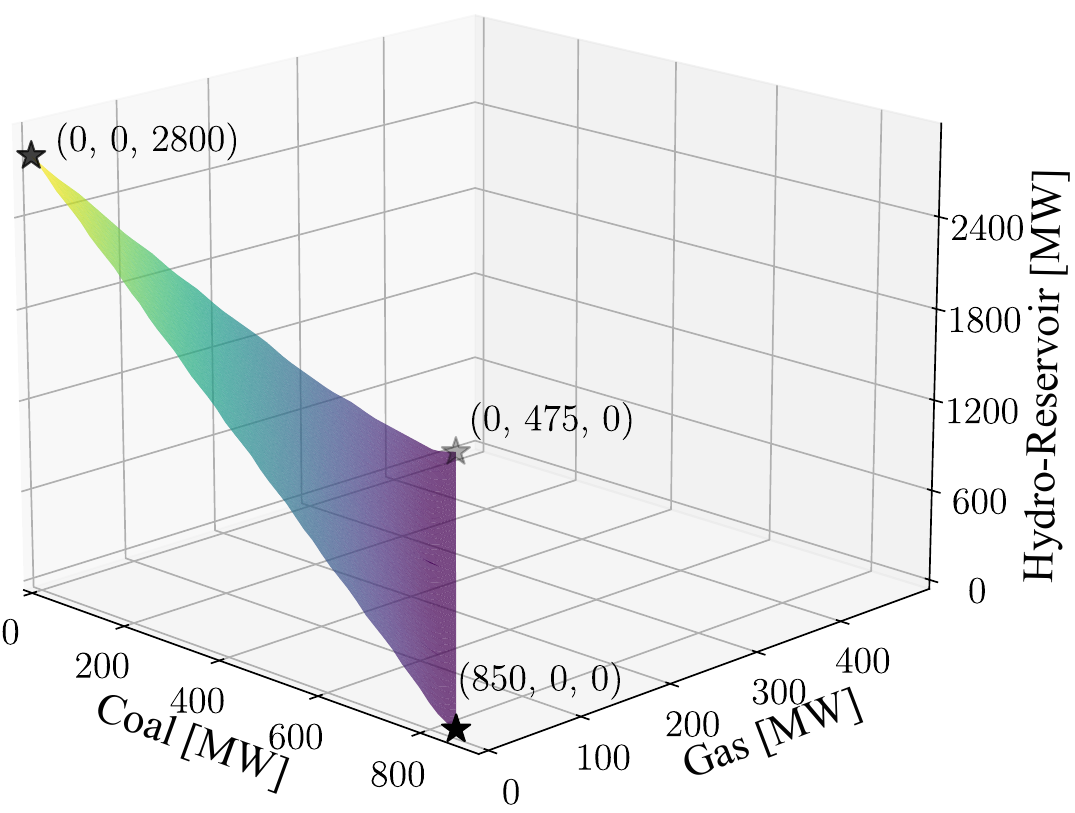}
    \caption{Nadir boundary for average demand in 2024}
    \label{fig:Figure_4}
\end{figure}

This process was carried out for the generation technologies that actively participated in frequency regulation in the Chilean electrical system in year 2024, namely coal-fired steam turbines, gas-fired combined-cycle turbines, and hydro-reservoirs. As an example, the boundary resulting for an average demand level during year 2024 is presented in Figure \ref{fig:Figure_4}, where it can be observed that the boundary can be conservatively approximated by a plane constructed from the edge points of the Nadir boundary (indicated by stars in the figure).

Then, for each hour in the FCUC problem where the Nadir requirement is not met, edge points are obtained through simulations, and a linear approximation of the boundary is derived and incorporated into the model as a Nadir constraint, as follows:

\begin{equation} \label{ec:nadir_methodology}
\begin{split}
    & b^{S}_t \sum_{i \in \mathcal{I}_\text{Coal}} \overline{P}_i u_{it} + b_t^{CC}\sum_{i \in \mathcal{I}_\text{Gas}} \overline{P}_i u_{it}  + b_t^{R}\sum_{h \in \mathcal{H}_\text{R}} \overline{P}_h u_{ht} - b^0_t \geq 0 
    \end{split}
\end{equation}

This approach is also applied for year 2035, which corresponds to a one-dimensional case, as only hydro-reservoirs are actively involved in frequency regulation in this case.

%% file: computational_experiments.tex
\subsection{General Setting and Case Studies} \label{sec:general_settings}

A single-node representation of the Chilean electrical system is developed, considering two contrasting operational contexts: years 2024 and 2035. The 2024 scenario reflects a system with a high share of conventional thermal units, whereas the 2035 scenario projects a near-complete retirement of these units, resulting in a grid dominated by IBRs and hydro-based generation \cite{CEN_2035}. Installed capacities for 2024 are based on real operational data \cite{CNE_2024, CEN_2024}, while those for 2035 are derived from the forescasts in \cite{CEN_2035}, adjusted to ensure system reliability under critical conditions. For both years, multiple scenarios are evaluated, defined according to the season (e.g., autumn or spring) and hydrological condition (e.g., dry or wet), with the aim of capturing relevant variations in renewable generation and water availability. Regarding the FCUC models studied, they are defined as follows:
\begin{itemize} 
    \item Industry Standard Model: Based on the formulation currently used in the Chilean industry \cite{CEN_métricas}, where a total spinning reserve requirement and inertia requirement are used. To satisfy the remaining frequency constraints (i.e., Nadir and QSS), the reserve requirement is iteratively increased until all frequency metrics fall within the threshold defined by the technical standards (see Table \ref{tab:lims_sen}).

    \input{Table_1}
    
    \item Proposed Model: Based on the proposed approach described in Section \ref{sec:methodology}, which uses a data-driven FCUC methodology and accurately embeds hydro-reservoir dynamics.
\end{itemize} 

To evaluate these models, a generation contingency of 400 MW is considered, consistent with the values employed by the Chilean system operator \cite{CEN_métricas}. The most relevant dynamic and operational parameters used are summarized in Table \ref{tab:dynamcs_params} and Table \ref{tab:operational_params}.
The models were implemented in the Julia programming language \cite{Julia} and the code is publicly available on GitHub \cite{code_github}. 
For the FCUC model, the JuMP \cite{JuMP} library and the Gurobi solver \cite{Gurobi} were used. Meanwhile, the dynamic model for the system frequency response (see Fig. \ref{fig:Figure_2}) was developed using the ModelingToolkit \cite{MTK} and DifferentialEquations \cite{Differentialequations} libraries.
Simulations were performed on a PC with an Intel Core i9 processor (2.60 GHz), 16 GB of RAM, and a 64-bit Windows 11 operating system.

\input{Table_2}

\input{Table_3}

The case studies conducted for this work are structured as follows. First, Section \ref{sec:dynamic_analysis} presents a dynamic analysis that examines the impact of various technologies on frequency constraints, with particular emphasis on the Nadir response.
Then, Section \ref{sec:operation_analysis} introduces an operational analysis, where the Chilean electrical system is simulated under different scenarios in which neither GFM inverters nor SCs are integrated, establishing a baseline under the current Chilean regulation \cite{norma_tecnica_CNE_seguridad_calidad}
Finally, considering ongoing technological developments, Section \ref{sec:planification_analysis} performs a planning analysis to evaluate the potential benefits of GFM inverters and SCs for frequency support.

\footnotetext{Represents battery degradation costs due to charge/discharge cycles, as per \cite{Hu}.}

\subsection{Case 1: Nadir Frequency Dynamic Analysis} \label{sec:dynamic_analysis}

The following section presents the dynamic analysis results for the generation mixes of years 2024 and 2035. These results were obtained by evaluating the dynamic model shown in Figure \ref{fig:Figure_2} for different values of the online capacity of the technologies considered in each year. 

For year 2024, the impact on the Nadir from reservoir hydropower, coal-fired units, and gas-fired units is analyzed, as these technologies had the largest share in frequency regulation during that year \cite{informe_desempeño}.
As shown in Figure \ref{fig:Figure_4}, the simulations reveal that combined-cycle gas turbines contribute the most significantly to enhancing the Nadir frequency.  
When analyzing the approximate relationships between technologies, it is observed that the effect on the Nadir of 1 MW of gas is equivalent to that of 1.8 MW of coal or 5.9 MW of hydro-reservoir. 

\begin{figure}[h!]
    \centering
    \includegraphics[width=0.65\linewidth]{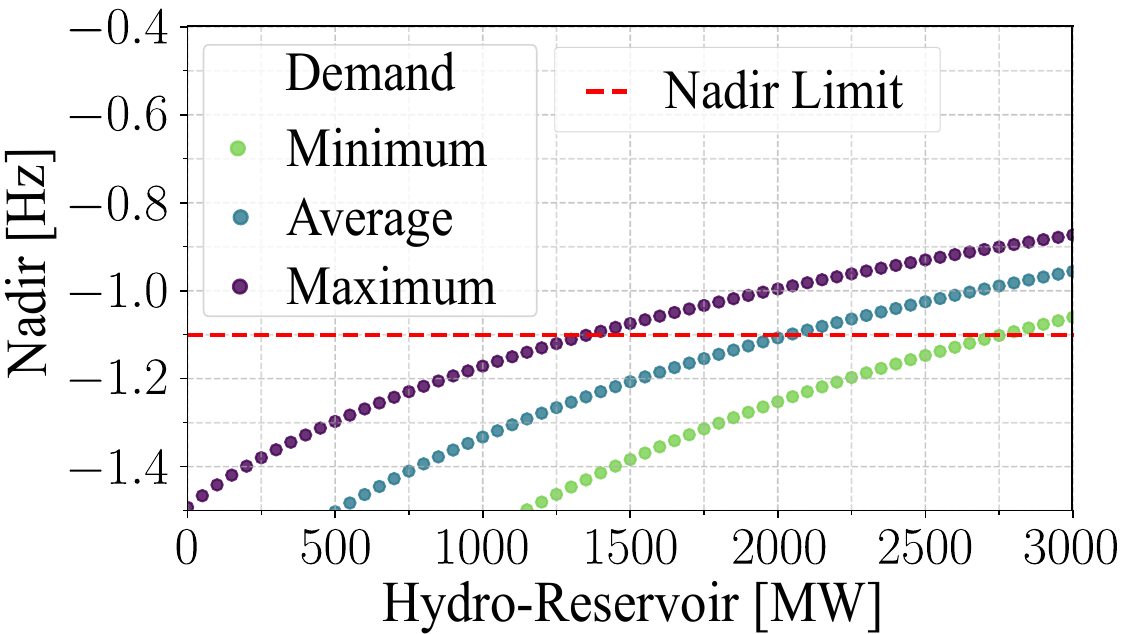}
    \caption{Impact of online capacity from hydro-reservoirs on Nadir frequency}
    \label{fig:Figure_5}
\end{figure}

For year 2035, given the retirement of thermal units and under the assumption that all inverters in the system are of the GFL type, hydro-reservoirs become the main technology contributing to frequency regulation. Based on this, Nadir frequency simulations are performed for different levels of demand and hydro-reservoir online capacity. The results are presented in Figure \ref{fig:Figure_5}. 

From the figure, it can be observed that increasing the online capacity of hydro-reservoirs improves the Nadir frequency. However, this relationship is not linear: the impact on Nadir diminishes as their online capacity in the system increases. On the other hand, higher demand levels reduce the minimum reservoir capacity required to meet the Nadir requirement, due to the increase in load damping’s contribution to frequency support.

\subsection{Case 2: Operation Models Analysis} \label{sec:operation_analysis}
The following section presents the results of the operational analysis comparing the industry standard and proposed models. For the 2024 scenario, Section \ref{sec:dry_autumn_2024} shows the results for autumn conditions with dry hydrology, while Section \ref{sec:wet_spring_2024} presents the results for spring conditions with wet hydrology. Similarly, for the 2035 scenario, Section \ref{sec:dry_autumn_2035} covers the results for dry autumn, and Section \ref{sec:wet_spring_2035} for wet spring. These conditions represent two significantly different scenarios in terms of renewable and water resource availability in the Chilean electrical system: dry autumn is characterized by low resource availability, whereas wet spring features a high resource availability.

\subsubsection{Dry Autumn 2024} \label{sec:dry_autumn_2024}
\begin{figure}[h]
    \centering
    \includegraphics[width=0.93\linewidth]{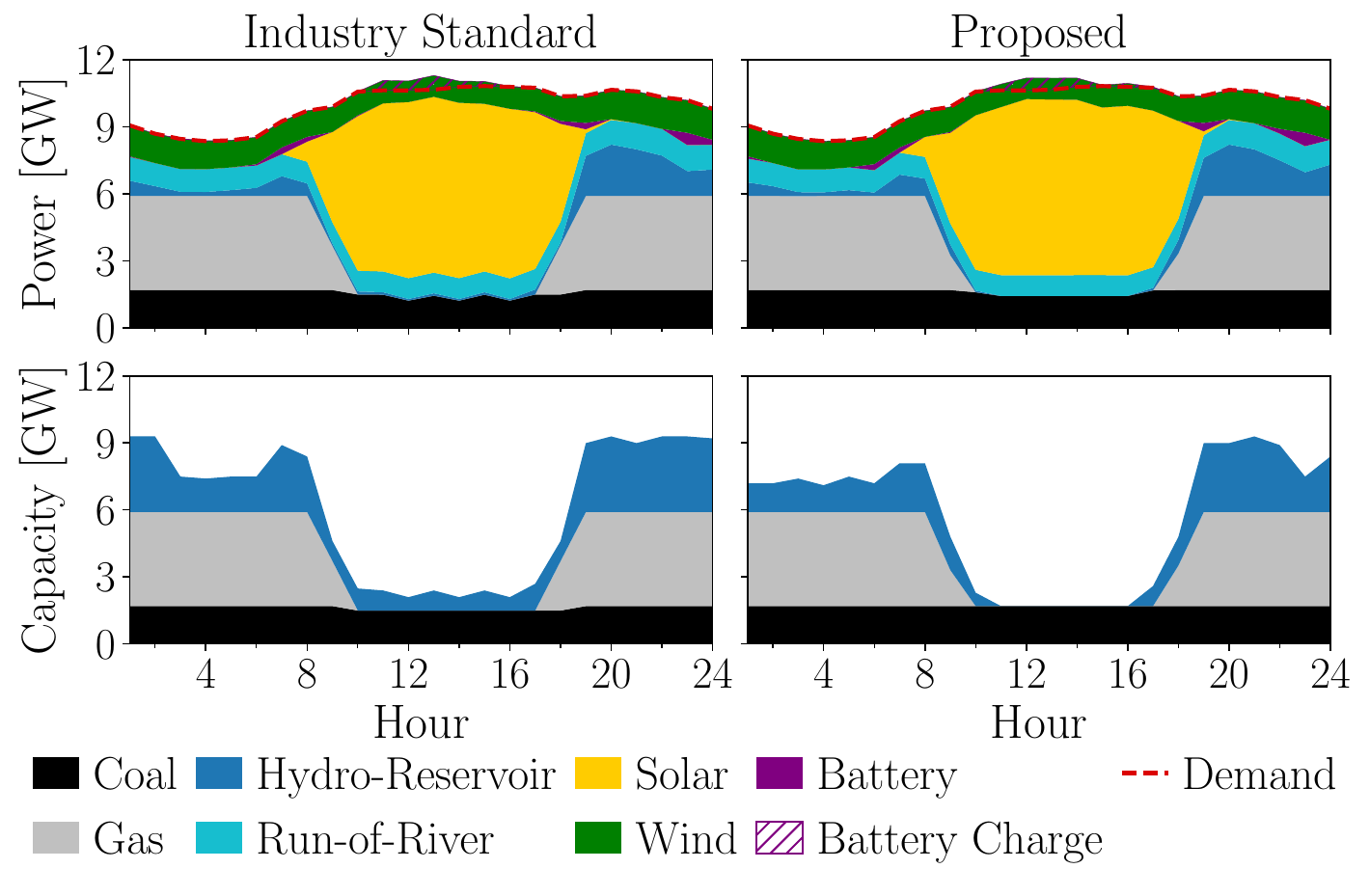}
    \caption{Power and capacity dispatch results for Dry-Autumn 2024 scenario}
    \label{fig:Figure_6}
\end{figure}

The results of the Dry Autumn 2024 scenario are presented in Table \ref{tab:costs_DA_2024} and Figure \ref{fig:Figure_6}. 
On the one hand, from Table \ref{tab:costs_DA_2024}, it can be observed that the proposed model achieves the lowest total daily operating cost, with only a slight difference of 0.11\% (13 kUSD) compared to the industry standard model. This is mainly associated with lower costs related to the operation of thermal units. 
\input{Table_4}

On the other hand, from Figure \ref{fig:Figure_6}, it can be seen that the proposed model is able to meet frequency requirements with a lower online capacity from hydro-reservoirs, relative to the industry standard model. This is because the industry standard approach increases reserve requirements uniformly throughout the day until frequency criteria are met in all hours, which leads to an overestimation of the number of online units needed. In contrast, the proposed model adjusts requirements on an hourly basis and incorporates the contribution of load damping in the system’s response, enabling a more efficient dispatch of generating units.

\subsubsection{Wet Spring 2024} \label{sec:wet_spring_2024}

\begin{figure}[h]
    \centering
    \includegraphics[width=0.93\linewidth]{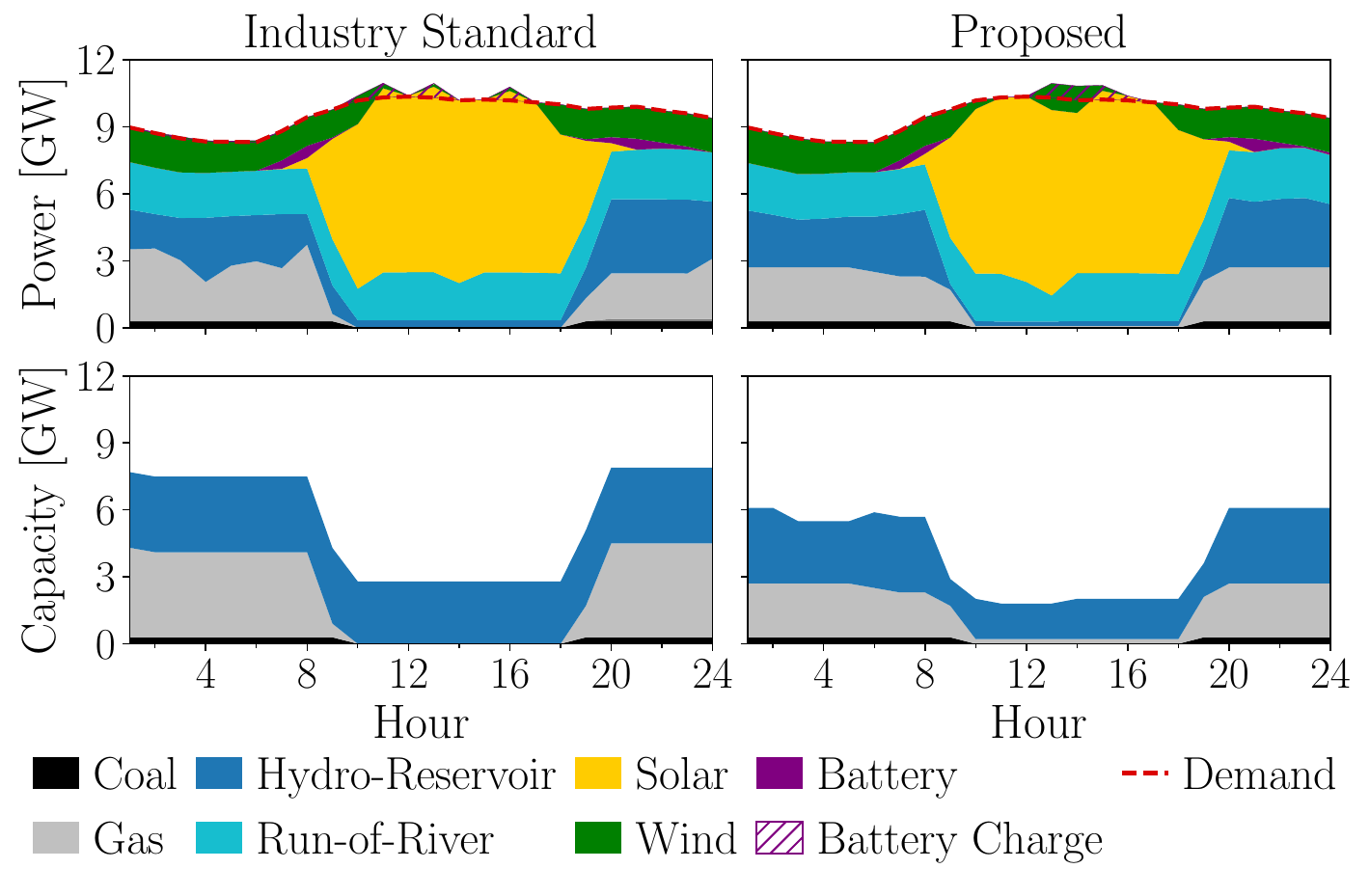}
    \caption{Power and capacity dispatch results for Wet-Spring 2024 scenario}
    \label{fig:Figure_7}
\end{figure}

The results of the Dry Autumn 2024 scenario are presented in Table \ref{tab:costs_WS_2024} and Figure \ref{fig:Figure_7}. 
On the one hand, from Table \ref{tab:costs_WS_2024}, it can be observed that the proposed model achieves the lowest total daily operating cost, with a 12\% difference (754 kUSD) compared to the industry standard model — a significantly larger gap than in the previous scenario. This is mainly associated with lower costs related to the operation of thermal units. 
\input{Table_5}

On the other hand, from Figure \ref{fig:Figure_7}, it can be seen that between hours 9 and 19, the proposed model dispatches fewer gas-fired thermal units compared to the industry standard model. As mentioned in the previous scenario, this is due to the industry standard model overestimating the number of online units required. In contrast, between hours 10 and 18—the most critical period due to the displacement of synchronous units by the high availability of renewable energy—the proposed model chooses to keep a gas unit operating at its minimum technical output to provide grid support. In comparison, the industry standard model meets the frequency requirements solely by dispatching hydro-reservoir units. This difference is explained by the fact that, during these hours, the proposed model introduces Nadir constraints that account for the specific impact of each technology on frequency regulation.  

These results demonstrate that the proposed model enables a more strategic use of generation technologies for secure system operation, adapting the dispatch to the specific requirements of each hour of the day and allowing for more cost-efficient outcomes compared to the industry standard model.

\subsubsection{Dry Autumn 2035} \label{sec:dry_autumn_2035}
\begin{figure}[h]
    \centering
    \includegraphics[width=0.93\linewidth]{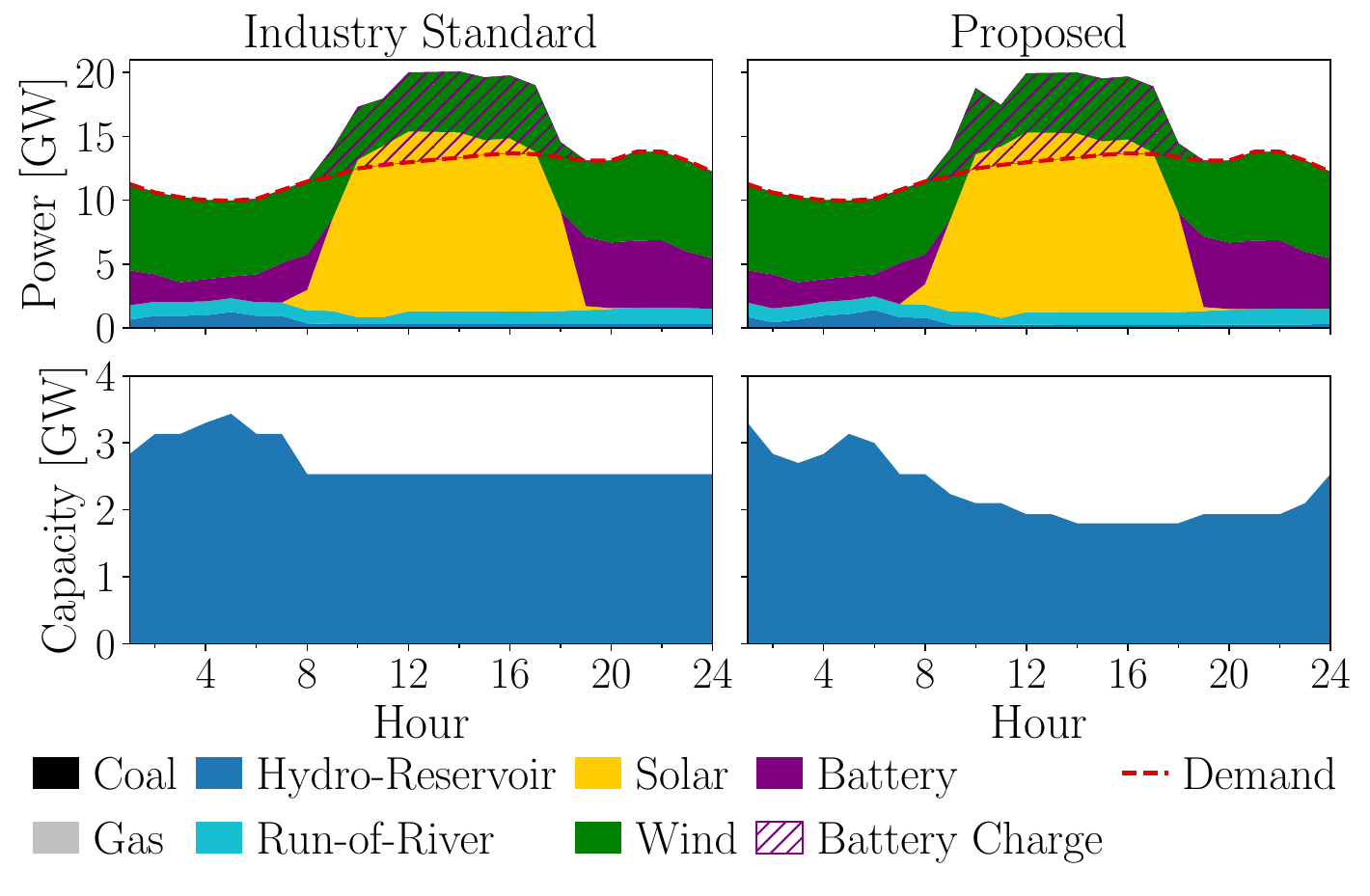}
    \caption{Power and capacity dispatch results for Dry-Autumn 2035 scenario}
    \label{fig:Figure_8}
\end{figure}

The results of the Dry Autumn 2024 scenario are presented in Table \ref{tab:costs_DA_2035} and Figure \ref{fig:Figure_8}. 
On the one hand, from Table \ref{tab:costs_DA_2035}, it can be observed that the proposed model achieves the lowest total daily operating cost, with a 3.6\% difference (96 kUSD) compared to the industry standard model. This is mainly associated with lower costs related to the operation of hydro-reservoirs. 
\input{Table_6}

On the other hand, from Figure \ref{fig:Figure_8}, it can be seen that the proposed model is able to meet frequency requirements with a lower online capacity from hydro-reservoirs compared to the industry standard model, similarly to what was observed in the dry-autumn 2024 scenario.

\subsubsection{Wet Spring 2035} \label{sec:wet_spring_2035}

\begin{figure}[h]
    \centering
    \includegraphics[width=0.93\linewidth]{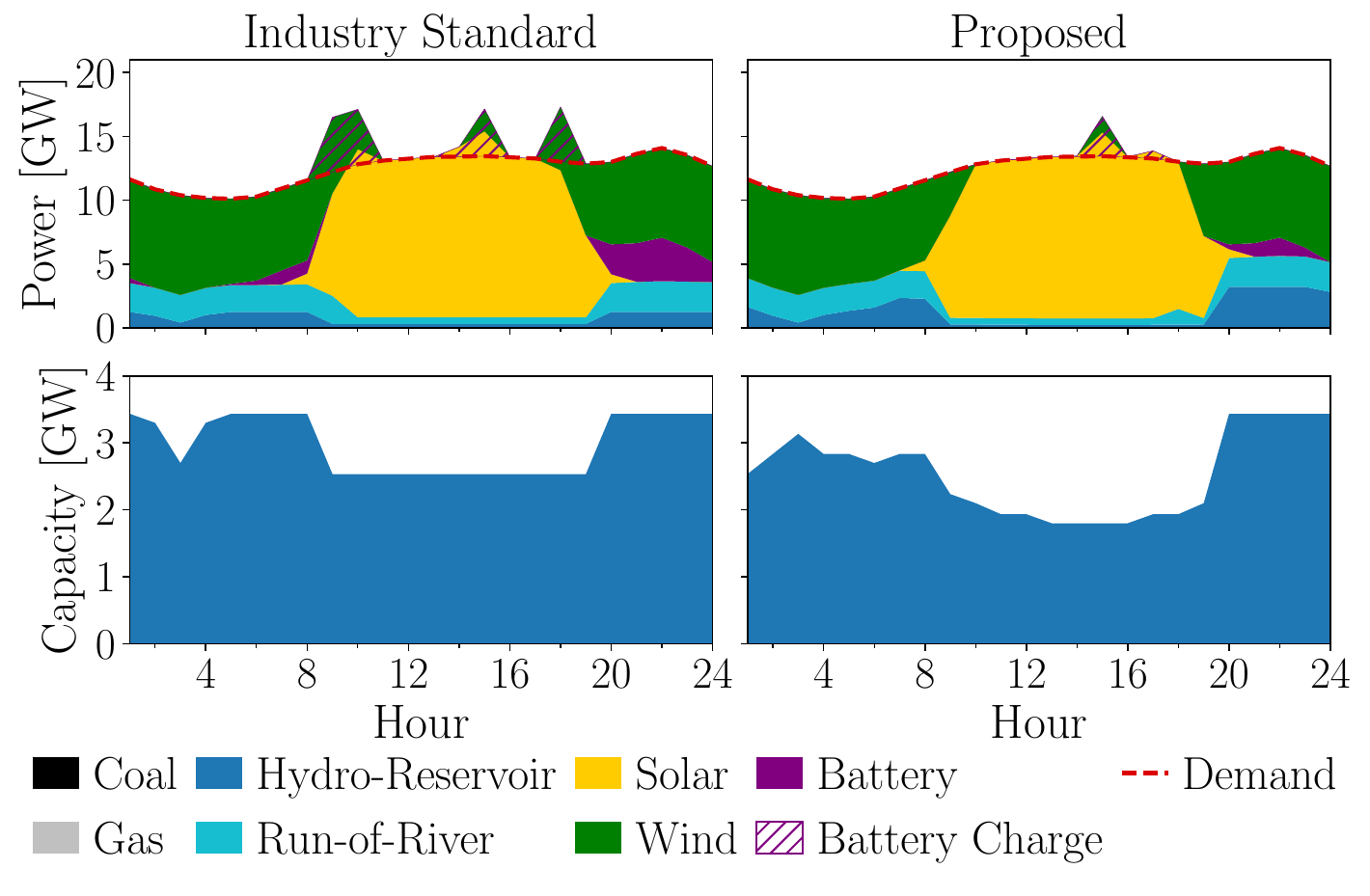}
    \caption{Power and capacity dispatch results for Wet-Spring 2035 scenario}
    \label{fig:Figure_9}
\end{figure}

The results of the Dry Autumn 2024 scenario are presented in Table \ref{tab:costs_WS_2035} and Figure \ref{fig:Figure_9}. 
On the one hand, from Table \ref{tab:costs_WS_2035}, it can be observed that the proposed model achieves the lowest total daily operating cost, with a 28\% difference (213 kUSD) compared to the industry standard model — a significantly larger gap than in the previous scenario. Although the proposed model has higher costs associated with reservoir hydropower generation, such cost is offset by a reduction in battery degradation costs resulting from charge/discharge cycles.
\input{Table_7}

On the other hand, from Figure \ref{fig:Figure_9}, it can be seen that the proposed model requires significantly less battery participation in system operation compared to the industry standard model.  
This contrast is primarily explained by the limitations of the industry model, which restricts the operational flexibility of reservoir hydropower units.
As a result, the system compensates with a more intensive use of batteries to meet frequency and energy balance requirements.
These results reaffirm that the proposed model enables a more strategic use of generation technologies compared to the industry standard approach.

\subsection{Case 3: Planning Analysis} \label{sec:planification_analysis}
The following section analyses the impact of installing grid-supporting technologies to improve frequency Nadir in 2035, namely in the form of Synchronous Condensers (SCs) and Grid-Forming (GFM) inverters. 

\begin{figure}[h]
    \centering
    \includegraphics[width=0.85\linewidth]{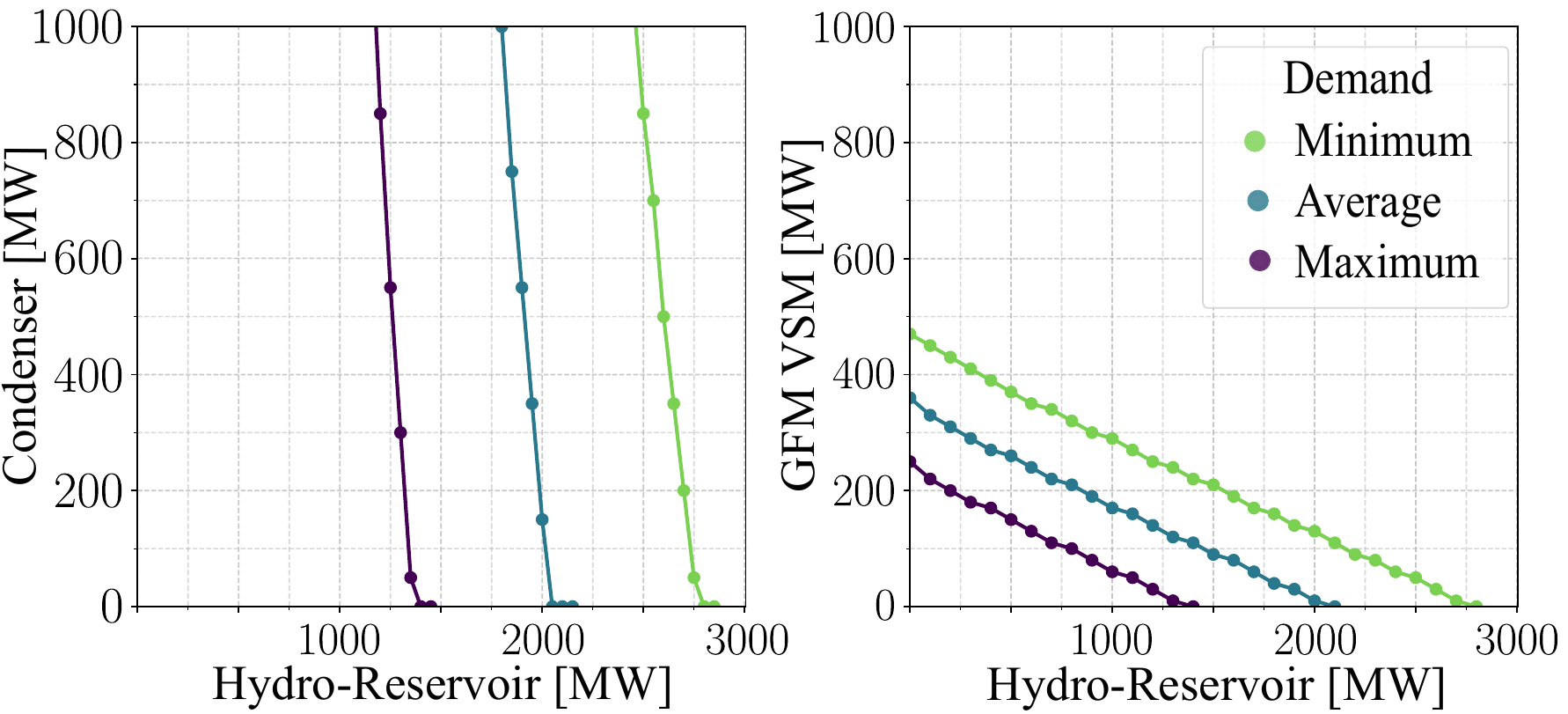}
    \caption{Nadir boundaries of SC and GFM VSM in relation to the spinning capacity of reservoir hydropower plants}
    \label{fig:Figure_10}
\end{figure}

First, a dynamic analysis was conducted to assess the impact of these grid-support technologies. As shown in Figure \ref{fig:Figure_10}, the results indicate that GFM inverters contribute most significantly to improving the Nadir frequency. When analyzing the approximate relationships between technologies, it is observed that the effect of 1 MW of GFM on Nadir is equivalent to that of 5.6 MW of hydro-reservoir or 30 MW of SC. 

 Second, a financial analysis is conducted by quantifying the annual operating cost savings achieved by installing either a 100 MW SC or a 100 MW GFM inverter of the VSM type. These savings are evaluated alongside investment costs associated with installation or retrofit, based on values reported in \cite{ISCI-ChileSustentable}. The assessment assumes a 6\% discount rate and a 20-year time horizon, and the results are summarized in Table~\ref{tab:financial_analysis}. The results show that the most attractive investment corresponds to GFM VSM inverters, particularly the retrofitting of existing units, which can be 28 times higher than SCs. Additionally, it is observed that, in the case of SCs, only retrofitting is financially viable, although it yields significantly lower profitability compared to the alternatives involving GFM VSM inverters. 
\input{Table_8}

\begin{figure}[h]
    \centering
    \includegraphics[width=0.45\linewidth]{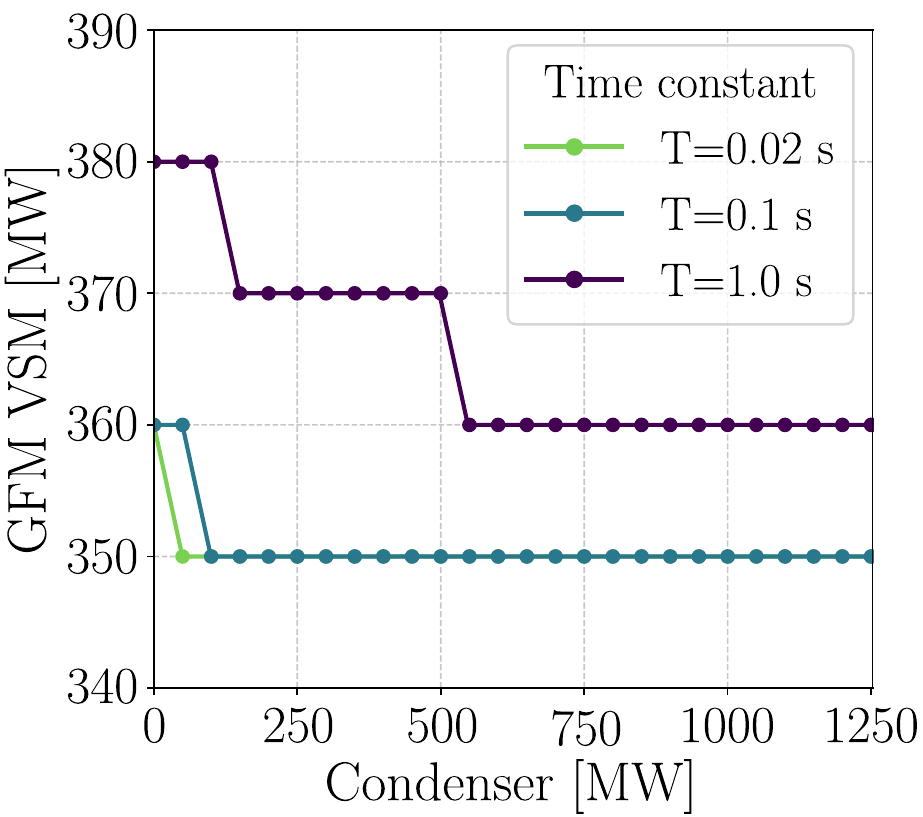}
    \caption{Nadir boundaries of SC and GFM VSM under different GFM time response constants}
    \label{fig:Figure_11}
\end{figure}

Finally, a sensitivity analysis of the response time of GFM inverters is made. From Table \ref{tab:dynamcs_params}, it can be observed that a time constant of 0.02 s was used in the simulations for GFM VSM inverters, as per [40]. However, other studies consider slower response times on the order of 0.1 s \cite{Pattabiraman, Nguyen}, and the Chilean system operator specifies a response within a range of up to 1 s \cite{CEN_2025_sensibilidad}. Therefore, the variation in the difference of the Nadir impact between GFM and SCs is evaluated for different inverter time constants. The resulting Nadir boundaries are shown in Figure \ref{fig:Figure_11}.

Based on these results, it is observed that for a time constant of 0.1 s, the requirement of 360 MW of GFM inverters to ensure compliance with the Nadir constraint—without the participation of SCs—remains unchanged, maintaining the equivalence that 1 MW of GFM has approximately the same impact on Nadir frequency as 30 MW of SCs. For a time constant of 1 s, the equivalence adjusts to 1 MW of GFM per 28 MW of SCs. Therefore, variations in the time constant do not significantly alter the simulation results.

%% file: Table_1.tex
\begin{table}[h]
\centering
\small
\renewcommand{\arraystretch}{1}
\begin{tabular}{lll}
    \hline
    \begin{tabular}[c]{@{}c@{}}Metric\end{tabular} & \begin{tabular}[c]{@{}c@{}}Threshold\end{tabular} & \begin{tabular}[c]{@{}c@{}}Reference Detail\end{tabular}  \\ \hline
    RoCoF  & 0.6 Hz/s & \cite{CEN_métricas} Section 1.1 \\ 
    Nadir  & 48.9 Hz  & \cite{CEN_métricas} Section 3.5 \\ 
    QSS    & 49.3 Hz  & \cite{CEN_métricas} Section 1.1 \\ \hline
\end{tabular}
\caption{Frequency metrics thresholds of the Chilean system operator.}
\label{tab:lims_sen}
\end{table}

%% file: Table_2.tex
\begin{table}[h]
\centering
\renewcommand{\arraystretch}{1}
\footnotesize
\begin{tabular}{llllllllllll}
\hline
\multicolumn{2}{l}{\begin{tabular}[l]{@{}l@{}}Steam\\ Turbine\end{tabular}} & \multicolumn{2}{l}{\begin{tabular}[l]{@{}l@{}}Combined-Cycle\\ Turbine\end{tabular}} & \multicolumn{2}{l}{\begin{tabular}[l]{@{}l@{}}Hydro\\ Reservoir\end{tabular}} & \multicolumn{2}{l}{\begin{tabular}[l]{@{}l@{}}GFM\\ VSM\end{tabular}} & \multicolumn{2}{l}{\begin{tabular}[l]{@{}l@{}}Run \\ of River\end{tabular}} & \multicolumn{2}{l}{\textbf{SC}} \\ \hline
M & 7.52 s & M & 13.14 s & M & 6.72 s & M & 1 s & M & 5.24 s & M & 10 s \\
R & 5\% & R & 5\% & R & 5\% & D & 20 & & & & \\
T & 7 s & T1 & 0.6 s & T & 1 s & T & 0.02 s & & & & \\
F & 0.3 & T2 & 0.5 s & Rtr & 38\% & & & & & & \\
& & T3 & 0.01 s & Ttr & 5 s & & & & & & \\
& & T4 & 0.24 s & & & & & & & & \\ \hline
\end{tabular}
\caption{Dynamic model parameters \cite{Kundur, FERNANDEZGUILLAMON, Unifi}}
\label{tab:dynamcs_params}
\end{table}

%% file: Table_3.tex
\begin{table}[h]
\centering
\renewcommand{\arraystretch}{1.1}
\footnotesize
\begin{tabular}{lllllllll}
\hline
Type & \begin{tabular}[l]{@{}l@{}}$P_{max}$\\ {[}MW{]}\end{tabular} & \begin{tabular}[l]{@{}l@{}}$P_{min}$\\ {[}MW{]}\end{tabular} & \begin{tabular}[l]{@{}l@{}}$C^V$\\ {[}kUSD{]}\end{tabular} & \begin{tabular}[l]{@{}l@{}}$C^{SU}$\\ {[}kUSD{]}\end{tabular} & \begin{tabular}[l]{@{}l@{}}$C^{SD}$\\ {[}USD{]}\end{tabular} & \begin{tabular}[l]{@{}l@{}}$C^F$\\ {[}USD{]}\end{tabular} & \begin{tabular}[l]{@{}l@{}}$T_{min}^{ON}$\\ {[}s{]}\end{tabular} & \begin{tabular}[l]{@{}l@{}}$T_{min}^{OFF}$\\ {[}s{]}\end{tabular} \\ \hline
Coal      & 200 & 80 & 0.06       & 60 & 4 & 12  & 8 & 8 \\
Gas       & 200 & 80 & 0.10       & 20 & 3 & 4   & 4 & 4 \\
Diesel    & 200 & 80 & 0.25       & 14 & 3 & 2.8 & 0 & 0 \\
Reservoir & 300 & 36 & 0.015-0.10 & -  & - & -   & - & - \\
Storage   & 100 & 0  & 0.015\footnotemark & -  & - & -   & - & - \\ \hline
\end{tabular}
\caption{Operational parameters \cite{PELP, PLP, ieee14bus, Córdova}}
\label{tab:operational_params}
\end{table}

%% file: Table_4.tex
\begin{table}[h]
\centering
\renewcommand{\arraystretch}{1}
\small
\begin{tabular}{lllllll}
\hline
\multirow{2}{*}{Model} & \multicolumn{4}{l}{Costs {[}kUSD{]}} & \multirow{2}{*}{\begin{tabular}[l]{@{}l@{}}Gap\\ {[}\%{]}\end{tabular}} & \multirow{2}{*}{\begin{tabular}[l]{@{}l@{}}Time\\ {[}min{]}\end{tabular}} \\ \cline{2-5} & 
\multicolumn{1}{l}{Thermal} & \multicolumn{1}{l}{Reservoir} & \multicolumn{1}{l}{Battery} & Total &  &  \\ \hline
\begin{tabular}[l]{@{}l@{}}Industry\\ Standard\end{tabular} & \multicolumn{1}{l}{12,591} & \multicolumn{1}{l}{1,497} & \multicolumn{1}{l}{59} & 14,147 & \multicolumn{1}{l}{0.1} & \multicolumn{1}{l}{17} \\ 
Proposed & \multicolumn{1}{l}{12,579} & \multicolumn{1}{l}{1,497} & \multicolumn{1}{l}{58} & 14,134 & \multicolumn{1}{l}{0.2} & \multicolumn{1}{l}{11} \\ \hline
\end{tabular}
\caption{Costs, optimality gap and simulation time for models in Dry-Autumn 2024 scenario}
\label{tab:costs_DA_2024}
\end{table}

%% file: Table_5.tex
\begin{table}[h]
\centering
\renewcommand{\arraystretch}{1}
\small
\begin{tabular}{lllllll}
\hline
\multirow{2}{*}{Model} & \multicolumn{4}{l}{Costs {[}kUSD{]}} & \multirow{2}{*}{\begin{tabular}[l]{@{}l@{}}Gap\\ {[}\%{]}\end{tabular}} & \multirow{2}{*}{\begin{tabular}[l]{@{}l@{}}Time\\ {[}min{]}\end{tabular}} \\ \cline{2-5} & 
\multicolumn{1}{l}{Thermal} & \multicolumn{1}{l}{Reservoir} & \multicolumn{1}{l}{Battery} & Total &  &  \\ \hline
\begin{tabular}[l]{@{}l@{}}Industry\\ Standard\end{tabular} &  \multicolumn{1}{l}{5,644} & \multicolumn{1}{l}{564} & \multicolumn{1}{l}{69} & 6,277 & 0.1 & 47 \\ 
Proposed & \multicolumn{1}{l}{4,898} & \multicolumn{1}{l}{564} & \multicolumn{1}{l}{61} & 5,523 & 0.2 & 4 \\ \hline
\end{tabular}
\caption{Costs, optimality gap and simulation time for models in Wet-Spring 2024 scenario}
\label{tab:costs_WS_2024}
\end{table}

%% file: Table_6.tex
\begin{table}[h]
\centering
\renewcommand{\arraystretch}{1}
\small
\begin{tabular}{lllllll}
\hline
\multirow{2}{*}{Model} & \multicolumn{4}{l}{Costs {[}kUSD{]}} & \multirow{2}{*}{\begin{tabular}[l]{@{}l@{}}Gap\\ {[}\%{]}\end{tabular}} & \multirow{2}{*}{\begin{tabular}[l]{@{}l@{}}Time\\ {[}min{]}\end{tabular}} \\ \cline{2-5} & 
\multicolumn{1}{l}{Thermal} & \multicolumn{1}{l}{Reservoir} & \multicolumn{1}{l}{Battery} & Total &  &  \\ \hline
\begin{tabular}[l]{@{}l@{}}Industry\\ Standard\end{tabular} & \multicolumn{1}{l}{0} & \multicolumn{1}{l}{1,179} & \multicolumn{1}{l}{1,487} & 2,666 & 0.2 & 31 \\ 
Proposed & \multicolumn{1}{l}{0} & \multicolumn{1}{l}{1,073} & \multicolumn{1}{l}{1,497} & 2,570 & 0.0 & 3 \\ \hline
\end{tabular}
\caption{Costs, optimality gap and simulation time for models in Dry-Autumn 2035 scenario}
\label{tab:costs_DA_2035}
\end{table}

%% file: Table_7.tex
\begin{table}[h]
\centering
\renewcommand{\arraystretch}{1}
\small
\begin{tabular}{lllllll}
\hline
\multirow{2}{*}{Model} & \multicolumn{4}{l}{Costs {[}kUSD{]}} & \multirow{2}{*}{\begin{tabular}[l]{@{}l@{}}Gap\\ {[}\%{]}\end{tabular}} & \multirow{2}{*}{\begin{tabular}[l]{@{}l@{}}Time\\ {[}min{]}\end{tabular}} \\ \cline{2-5} & 
\multicolumn{1}{l}{Thermal} & \multicolumn{1}{l}{Reservoir} & \multicolumn{1}{l}{Battery} & Total &  &  \\ \hline
\begin{tabular}[l]{@{}l@{}}Industry\\ Standard\end{tabular} &  \multicolumn{1}{l}{0} & \multicolumn{1}{l}{270} & \multicolumn{1}{l}{500} & 770 & 0.1 & 0.5 \\ 
Proposed & \multicolumn{1}{l}{0} & \multicolumn{1}{l}{445} & \multicolumn{1}{l}{112} & 557 & 0.0 & 2 \\ \hline
\end{tabular}
\caption{Costs, optimality gap and simulation time for models in Wet-Spring 2035 scenario}
\label{tab:costs_WS_2035}
\end{table}

%% file: Table_8.tex
\begin{table}[h]
\centering
\renewcommand{\arraystretch}{1.3}
\footnotesize
\begin{tabular}{lllll}
\cline{2-5}
\multicolumn{1}{l}{} & \multicolumn{2}{l}{SC} & \multicolumn{2}{l}{GFM VSM} \\ \cline{2-5} & 
\multicolumn{1}{l}{Retrofit} & New unit & \multicolumn{1}{l}{Retrofit} & New unit \\ \hline
\multicolumn{1}{l}{VAN {[}kUSD/MW{]}} & \multicolumn{1}{l}{45} & -264 & \multicolumn{1}{l}{1,266} & 1,211 \\ 
\multicolumn{1}{l}{TIR {[}\%{]}} & \multicolumn{1}{l}{12} & -5 & \multicolumn{1}{l}{115} & 75 \\ 
\multicolumn{1}{l}{PayBack {[}years{]}} & \multicolumn{1}{l}{7.7} & 33.6 & \multicolumn{1}{l}{0.9} & 1.3 \\ \hline
\end{tabular}
\caption{Financial results of GFM inverter and SC integration analysis}
\label{tab:financial_analysis}
\end{table}

%% file: conclusion.tex
This paper has addressed the need to incorporate frequency requirements into the daily operation problem in power systems with significant hydro participation. To this end, a dynamic model was developed to quantify the impact of different generation and/or grid-support technologies on meeting frequency requirements, including hydro-reservoirs and GFM-based plants. Additionally, an FCUC model that accounts for the dynamic characteristics of the different system units was formulated and evaluated for the Chilean grid during years 2024 (current scenario) and 2035 (carbon-free scenario). 

From a strictly dynamic point of view, combined-cycle gas units show overall better performance compared to other conventional technologies based on synchronous machines. However, when the economic perspective is added to the analysis, it is revealed that hydro-reservoirs are a fundamental technology due to their relatively low operational cost when compared to gas units. These results highlight the relevance of hydro-reservoirs in frequency regulation, with these units expected to play a critical role as thermal units are retired in the coming years. The proposed model demonstrated that incorporating frequency constraints while accurately modeling the dynamic response of hydro-based generation units can reduce daily operational costs by up to 28\% for the Wet Autumn 2024 scenario compared to the industry standard model. These cost reductions are more pronounced in case studies with high participation of IBRs, driven by greater availability of renewable generation. 

Furthermore, a planning analysis was carried out to assess the economic feasibility of integrating GFM inverters and SCs to support frequency regulation. The results indicated that GFM significantly outperforms SCs in terms of investment return. In the case of retrofitting existing units, GFMs provide 28 times higher profitability, and for new installations, only GFMs achieve a positive return. 

Based on these findings, the following recommendations are proposed to strengthen the operation of the Chilean electrical system in the coming years: 
\begin{itemize}
    \item Frequency constraints should be incorporated into the planning problem for the use of reservoir-based hydro resources, given the critical role these units will play following the retirement of thermal power plants.
    \item Technologies capable of supporting frequency regulation should be assessed to avoid overreliance on hydro-reservoirs. A potential alternative is to keep combined cycle gas units available, given their high impact on meeting frequency requirements. Another option is the integration of GFM VSM inverters, due to both their efficiency in frequency regulation and the economic feasibility of their development.
\end{itemize}

For future work, it is proposed to develop a hydroelectric planning problem that explicitly incorporates frequency constraints. Additionally, a multinodal case study is suggested, allowing the dynamic response to be represented in a disaggregated manner, thus enabling an assessment of frequency requirements at different system areas. This is particularly relevant for systems exhibiting a string-like topology, as in the case of Chile, which is non-meshed and primarily longitudinal \cite{CEN2025b}.

%% file: CRediT.tex
\textbf{Valeria Aravena:} Conceptualization, Methodology, Software, Writing – Original Draft.
\textbf{Samuel Córdova:} Supervision, Conceptualization, Methodology, Writing - Review \& Editing.
\textbf{Maximiliano Kairath:} Methodology, Software, Review \& Editing.
\textbf{Matías Negrete-Pincetic:} Supervision, Review \& Editing.